\documentclass[11pt]{llncs}

\usepackage{diagrams}
\usepackage{fancybox}
\usepackage{amsmath}
\usepackage{amssymb}
\usepackage{amsfonts}
\DeclareMathOperator{\charpoly}{charpoly}
\DeclareMathOperator{\splittemp}{split}
\renewcommand{\split}{\splittemp} \DeclareMathOperator{\deriv}{deriv}

\newcommand{\hra}{\hookrightarrow}

\DeclareMathOperator{\disc}{disc}
\def\Zbar{\overline{\Z}}
\def\GL{\mathrm{GL}}

\DeclareMathOperator{\tr}{tr}
\DeclareMathOperator{\SL}{SL}
\DeclareMathOperator{\End}{End}
\DeclareMathOperator{\new}{new}
\def\Tn{\T^{\new}}
\DeclareMathOperator{\ord}{ord}
\DeclareMathOperator{\Gal}{Gal}
\newcommand{\T}{\mathbf{T}}
\newcommand{\Z}{\mathbf{Z}}
\newcommand{\C}{\mathbf{C}}
\newcommand{\Q}{\mathbf{Q}}
\newcommand{\Qbar}{\overline{\Q}}

\newcommand{\Qpbar}{\overline{\Q}_p}
\newcommand{\tensor}{\otimes}
\newcommand{\tT}{\widetilde{\T}}
\newcommand{\tTn}{\tT^{\new}}
\newcommand{\F}{\mathbf{F}}
\newcommand{\Fbar}{\overline{\F}}
\newcommand{\Fpbar}{\overline{\F}_p}
\newcommand{\comment}[1]{}
\newcommand{\con}{\equiv}
\newcommand{\ds}{\displaystyle}
\def\Tw{\widetilde{\T}}
\def\rhobar{\overline{\rho}}

\title{Conjectures About Discriminants of Hecke Algebras of Prime Level}
\titlerunning{Conjectures About Discriminants}
\author{Frank Calegari\footnote{Supported in part by the American Institute
of Mathematics}\inst{1} \and William A. Stein\footnote{Supported in part
by a National Science Foundation Postdoctoral Fellowship}\inst{2}}
\authorrunning{Calegari \and Stein}
\tocauthor{Frank Calegari (Harvard University),
William Stein (Harvard University)}
\institute{Harvard University\\
\email{fcale@math.harvard.edu}\\
\texttt{http://www.math.harvard.edu/\~{}fcale}
\and
Harvard University,\\
\email{was@math.harvard.edu}\\
\texttt{http://modular.fas.harvard.edu/}}
\begin{document}
\maketitle

\begin{abstract}
  In this paper, we study $p$-divisibility of discriminants of Hecke
  algebras associated to spaces of cusp forms of prime level.  By
  considering cusp forms of weight bigger than~$2$, we are are led to
  make a precise conjecture about indexes of Hecke algebras in their
  normalisation which implies (if true) the surprising conjecture that
  there are no mod~$p$ congruences between non-conjugate newforms in
  $S_2(\Gamma_0(p))$, but there are almost always many such
  congruences when the weight is bigger than~$2$.
\end{abstract}

\section{Basic Definitions}
We first recall some commutative algebra related to discriminants,
then introduce Hecke algebras of spaces of cusp forms.

\subsection{Commutative Algebra}
In this section we recall the definition of discriminant of a finite
algebra and note that the discriminant is nonzero if and only if
no base extension of the algebra contains nilpotents.

Let~$R$ be a ring and let~$A$ be an~$R$-algebra that is free of finite
rank as an~$R$-module.  The {\em trace} of $x\in A$ is the trace, in
the sense of linear algebra, of left multiplication by~$x$.

\begin{definition}[Discriminant]
  Let $\omega_1,\ldots,\omega_n$ be an~$R$-basis for~$A$.  Then the
  {\em discriminant} $\disc(A)$ of~$A$ is the determinant of
  the $n\times n$ matrix $(\tr(\omega_i\omega_j))$. 
\end{definition}
The discriminant is only well-defined modulo squares of units in~$R$.
When $R=\Z$ the discriminant is well defined, since the only units are
$\pm 1$.

We say that $A$ is \textit{separable over $R$} if 
for every extension $R'$ of $R$, the 
ring $A\tensor R'$ contains no nilpotents.
\begin{proposition}\label{prop:separable}
  Suppose~$R$ is a field.  Then~$A$ has nonzero 
  discriminant if and only if~$A$ is 
  separable over~$R$.
\end{proposition}
\begin{proof}
  For the convenience of the reader, we summarize the proof in
  \cite[\S26]{matsumura}.  If~$A$ contains a nilpotent then that
  nilpotent is in the kernel of the trace pairing, so the discriminant
  is~$0$.  Conversely, if~$A$ is separable then we may assume that~$R$
  is algebraically closed.  Then~$A$ is an Artinian reduced ring,
  hence isomorphic as a ring to a finite product of copies of~$R$,
  since~$R$ is algebraically closed.  Thus the trace form on~$A$ is
  nondegenerate.
\end{proof}

\subsection{The Discriminant Valuation}
We next introduce Hecke algebras attached to certain spaces of cusp
forms of prime level $p$, define the discriminant valuation as the
exponent of the largest power of~$p$ that divides the discriminant,
and observe that there are eigenform congruences modulo~$p$ exactly
when the discriminant valuation is positive. We then present an
example to illustrate the definitions.

Let $\Gamma$ be a congruence subgroup of $\SL_2(\Z)$.  In this paper,
we will only consider $\Gamma=\Gamma_0(p)$ for~$p$
prime.  For any positive integer~$k$, let $S_k(\Gamma)$ denote the space
of holomorphic weight~$k$ cusp forms for~$\Gamma$.  Let
$$
\T = \Z[\ldots,T_n,\ldots] \subset \End(S_k(\Gamma))
$$
be the associated Hecke algebra, which is generated by Hecke
operators $T_n$ for all integers~$n$, including $n=p$ 
(we will sometimes write $U_p$ for $T_p$).
Then~$\T$ is a commutative ring that is free as a module over~$\Z$ of
rank equal to $\dim S_k(\Gamma)$.  We will also sometimes consider the
image $\T^{\new}$ of~$\T$ in $\End(S_k(\Gamma)^{\new})$.

\comment{
\begin{example}
  Let $\Gamma=\Gamma_0(243)$.
  Since $243=3^5$, there are nilpotents in $\T$, so $\disc(\T) =0$.  
  A computation shows that
  $$
  \disc(\T^{\new}) = 2^{13} \cdot 3^{40}.
  $$
\end{example}
}

\begin{definition}[Discriminant Valuation]
  Let~$p$ be a prime,~$k$ a positive integer, and suppose that $\Gamma=\Gamma_0(
p)$.
  Let $\T$ be the corresponding Hecke algebra.
Then the {\em discriminant valuation of $\Gamma$ in weight~$k$} is
  $$
  d_k(\Gamma) = \ord_p(\disc(\T)).
  $$
\end{definition}

We expect that $d_k(\Gamma)$ is finite for the following reason.
The Hecke operators $T_n$, with~$n$ not divisible by~$p$, are
diagonalizable since they are self adjoint with respect to the
Petersson inner product.  When $k=2$ one knows that $U_p$ is
diagonalizable since the level is square free, and when $k>2$ one
expects this (see \cite{coleman-edixhoven}).  If~$\T$ contains no
nilpotents, Proposition~\ref{prop:separable} implies that the
discriminant of $\T$ is nonzero.  Thus $d_k(\Gamma)$ is finite when
$k=2$ and conjectured to be finite when $k>2$.

Let~$p$ be a prime and suppose that $\Gamma=\Gamma_0(p)$.
 A {\em normalised eigenform} is an element $f=\sum a_n
q^n\in S_k(\Gamma)$ that is an eigenvector for {\em all} Hecke
operators $T_\ell$, including those that divide~$p$, normalised so
that $a_1 = 1$.  The quantity $d_k(\Gamma)$ is of interest because it
measures mod~$p$ congruences between normalised eigenforms in
$S_k(\Gamma)$.

\begin{proposition}
Assume that $d_k(\Gamma)$ is finite.  The discriminant
  valuation $d_k(\Gamma)$ is positive (i.e., the discriminant is
  divisible by~$p$) if and only if there is a congruence in
  characteristic~$p$ between two normalized eigenforms in
  $S_k(\Gamma)$. (The two congruent eigenforms might be Galois
  conjugate.)
\end{proposition}
\begin{proof}
  It follows from Proposition~\ref{prop:separable} that
  $d_k(\Gamma)>0$ if and only if $\T\tensor \Fpbar$ is not separable.
  The Artinian ring $\T\tensor\Fpbar$ is not separable if and only if
  the number of ring homomorphisms $\T\tensor\Fpbar \to \Fpbar$ is
  less than
  $$
  \dim_{\Fpbar} \T\tensor\Fpbar = \dim_\C S_k(\Gamma).
  $$
  Since $d_k(\Gamma)$ is finite, the number of ring homomorphisms
  $\T\tensor\Qpbar \to \Qpbar$ equals $\dim_\C S_k(\Gamma)$.  The
  proposition follows from the fact that for any ring~$R$, there is a
  bijection between ring homomorphisms $\T\to R$ and normalised
  eigenforms with $q$-expansion in~$R$.
\end{proof}
The same proof also shows that a prime $\ell$ divides the discriminant
of $\T$ if and only if there is a congruence mod~$\ell$ between two
normalized eigenforms in $S_k(\Gamma)$

\begin{example}
  If $\Gamma=\Gamma_0(389)$ and $k=2$, then $\dim_\C S_2(\Gamma) =
  32$.  Let~$f$ be the characteristic polynomial of $T_2$.  One can
  check that~$f$ is square free and $389$ exactly divides the
  discriminant of~$f$. This implies that $d_2(\Gamma)=1$ and that
  $T_2$ generates $\T\tensor \Z_{389}$ as an algebra over
  $\Z_{389}$. (If $T_2$ only generated a subring of
  $\T\tensor\Z_{389}$ of finite index $>1$, then the discriminant
  of~$f$ would be divisible by $389^2$.)
  
  Modulo~$389$ the characteristic polynomial~$f$ is congruent to
  $$\begin{array}{l}
    (x+2)(x+56)(x+135)(x+158)(x+175)^2(x+315)(x+342)(x^2+387)\\
    (x^2+97x+164)(x^2 + 231x + 64)(x^2 + 286x + 63)(x^5 + 88x^4 +196x^3 + \\
    113x^2 +168x + 349)(x^{11} + 276x^{10} + 182x^9 + 13x^8 + 298x^7 + 316x^6 +\\
    213x^5 + 248x^4 + 108x^3 + 283x^2 + x + 101)
  \end{array}
  $$
  The factor $(x+175)^2$ indicates that $\T\tensor \F_{389}$ is not
  separable over $\F_{389}$ since the image of
  $(\overline{f}/(x+175))(T_2)$ in $\T\tensor \F_{389}$ is nilpotent
  (it is nonzero but its square is~$0$).  There are $32$ eigenforms
  over~$\Q_2$ but only $31$ mod~$389$ eigenforms, so there must be a
  congruence.  There is a newform~$F$ in
  $S_2(\Gamma_0(389),\Zbar_{389})$ whose $a_2$ term is a root of
  \[
  x^2 + (-39 + 190\cdot 389 + 96\cdot 389^2 +\cdots) x + (-106 +
  43\cdot 389 + 19\cdot 389^2 + \cdots).
  \] 
  There is a congruence between~$F$ and its
  $\Gal(\Qbar_{389}/\Q_{389})$-conjugate.
\end{example}

\section{Computing Discriminants}
In this section we sketch the algorithm that we use for computing the
discriminants mentioned in this paper.  

This algorithm was inspired by a discussion of the second author with
Hendrik Lenstra.  We leave the details of converting the description
below into standard matrix operations to the reader.  Also, the
modular symbols algorithms needed to compute Hecke operators are quite
involved.

Let $\Gamma=\Gamma_0(p)$, and let $k\geq 2$ be an
integer.  The following sketches an algorithm for computing the
discriminant of the Hecke algebra $\T$ acting on $S_k(\Gamma)$.
\begin{enumerate}
\item For any given $n$, we can explicitly compute a matrix that
represents the action of Hecke operators $T_n$ on $S_k(\Gamma)$ using
modular symbols.  We use the second author's MAGMA \cite{magma}
packages for computing with modular symbols, which builds on work of
many people (including \cite{cremona:algs} and \cite{merel:1585}).
\item Using the Sturm bound, as described in the appendix to
\cite{agashe-stein:schoof-appendix}, find an integer~$b$ such that
$T_1,\ldots,T_b$ generate $\T$ as a $\Z$-module.  (The integer $b$
is $\lceil{}(k/12)\cdot [\SL_2(\Z):\Gamma]\rceil{}$.)

\item  Find a subset~$B$ of the $T_i$ that form a $\Q$-basis for
$\T\tensor_\Z\Q$.    (This uses Gauss elimination.)

\item View $\T$ as a ring of matrices acting on $\Q^{d}$, where
$d=\dim(S_k(\Gamma))$ and try random sparse vectors $v\in\Q^{d}$
until we find one 
such that the set of vectors $C=\{T(v) : T \in B\}$ are linearly
independent.

\item Write each of $T_1(v),\ldots, T_b(v)$ as $\Q$-linear combinations of the
elements of~$C$. 

\item Find a $\Z$-basis~$D$ for the $\Z$-span of
these $\Q$-linear combinations of elements of~$C$.  
(This basis $D$ corresponds to a $\Z$-basis for $\T$, but
is much easier to find that directly looking for
a $\Z$-basis in the space of $d\times d$ matrices
that $\T$ is naturally computed in.)

\item Unwinding what we have done in the previous steps, find the
  trace pairing on the elements of~$D$, and deduce the discriminant
  of~$\T$ by computing the determinant of the trace pairing matrix.
\end{enumerate}

A very time-consuming step, at least in our implementation, is
computing~$D$ from $T_1(v),\ldots,T_b(v)$ expressed in terms of~$C$,
and this explains why we embed $\T$ in $\Q^{d}$ instead of viewing the
elements of $\T$ as vectors in $\Q^{d \times d}$.

An implementation by the second author of the above algorithm is
included with the MAGMA computer algebra system.  The relevant source
code is in the file {\tt Geometry/ModSym/linalg.m} in the {\tt
package} directory (or ask the second author of the apper to send you
a copy {\tt linalg.m}).  We illustrate the use of MAGMA to compute
discriminants below, which were run under MAGMA V2.10-21 for Linux on
a computer with an Athlon 2800MP processor (2.1Ghz).
\begin{verbatim}
  > M := ModularSymbols(389,2, +1);
  > S := CuspidalSubspace(M);
  > time D  := DiscriminantOfHeckeAlgebra(S);
  Time: 0.750
  > D;
  629670054720061882880174736321392595498204931550235108311\
  04000000
  > Factorisation(D);
  [ <2, 53>, <3, 4>, <5, 6>, <31, 2>, <37, 1>, <389, 1>, ...]
  > M := ModularSymbols(997,2, +1); S := CuspidalSubspace(M);
  > time D  := DiscriminantOfHeckeAlgebra(S);
  Time: 55.600
\end{verbatim}
The reason for the $+1$ in the construction of modular symbols is so
that we compute on a space that is isomorphic as a $\T$-module to one
copy of $S_2(\Gamma_0(p))$, instead of two copies.

\section{Data About Discriminant Valuations}
In this section we report on our extensive computations of
$d_k(\Gamma_0(p))$.  We first note that there is only one~$p<50000$
such that $d_2(\Gamma_0(p))>0$.  Next we give a table of values of
$d_4(\Gamma_0(p))$, which seems to exhibit a nice pattern.

\subsection{Weight Two}
\begin{theorem}\label{thm:disc}
  The only prime $p<60000$ such that $d_2(\Gamma_0(p))>0$ is $p=389$,
  with the possible exception of $50923$ and $51437$.
\end{theorem}
Computations in this direction by the second author have been cited in
\cite{ribet:torsion}, \cite{merel-stein}, \cite{ono-mcgraw}, and
\cite{momose-ozawa}.  For example, Theorem~\ref{thm:disc} is used for
$p<1000$ in \cite{merel-stein} as a crucial step in proving that
if~$E$ is an elliptic curve over $\Q(\mu_p)$, with $17\leq p <1000$,
then not all elements of $E(\Qbar)[p]$ are rational over $\Q(\mu_p)$.

\begin{proof}
  This is the result of a large computer computation.  The rest of
  this proof describes how we did the computation, so the reader has
  some idea how to replicate or extend the computation.  The
  computation described below took about one week using a cluster
  equipped with $10$ Athlon 2000MP processors.  The computations are
  nontrivial; we compute spaces of modular symbols, supersingular
  points, and Hecke operators on spaces of dimensions up to~$5000$.

  
  The aim is to determine whether or not~$p$ divides the discriminant
  of the Hecke algebra of level~$p$ for each $p < 60000$.  If~$T$ is
  an operator with integral characteristic polynomial, we write
  $\disc(T)$ for $\disc(\charpoly(T))$, which also equals
  $\disc(\Z[T])$. We will often use that
  $$\disc(T)\!\!\!\!\mod{p} = \disc(\charpoly(T)\!\!\!\!\mod p).$$
  
  We ruled out the possibility that $d_k(\Gamma_0(p))>0$ for most
  levels~$p<60000$ by computing characteristic polynomials of Hecke
  operators using an algorithm that the second author and D.~Kohel
  implemented in MAGMA (\cite{magma}), which is based on the
  Mestre-Oesterle method of graphs \cite{mestre:graphs} (or contact
  the second author for an English translation).  Our implementation
  is available as the ``Module of Supersingular Points'' package that
  comes with MAGMA.  We computed $\disc(T_q)$ modulo~$p$ for several
  small primes~$q$, and in most cases found a prime~$q$ such that this
  discriminant is nonzero.  The following table summarises how often
  we used each prime~$q$ (note that there are $6057$ primes up to
  $60000$):
\begin{center}
\begin{tabular}{|l|l|}\hline
$q$  & number of $p< 60000$ where~$q$ smallest
  s.t. $\disc(T_q)\neq 0$ mod~$p$\\\hline
2&             5809 times\\
3&             161   (largest: 59471)\\
5&             43    (largest: 57793)\\
7&             15    (largest: 58699)\\
11&            15    (the smallest is 307; the largest 50971)\\
13&            2     (they are 577 and 5417)\\
17&            3     (they are 17209, 24533, and 47387)\\
19&            1     (it is 15661 )\\\hline
\end{tabular}
\end{center}

The numbers in the right column sum to 6049, so 8 levels are missing.
These are
$$
389,487,2341,7057,15641,28279, 50923, \text{ and } 51437.
$$
(The last two are still being processed.  $51437$ has the property
that $\disc(T_q)=0$ for $q=2,3,\ldots,17$.)  We determined the
situation with the remaining 6 levels using Hecke operators $T_n$
with~$n$ composite.
\begin{center}
\begin{tabular}{|l|l|}\hline
$p$ & How we rule level~$p$ out, if possible\\\hline
389&   $p$ does divide discriminant\\
487&   using charpoly($T_{12}$)\\
2341&  using charpoly($T_6$)\\
7057&  using charpoly($T_{18}$)\\
15641& using charpoly($T_6$)\\
28279& using charpoly($T_{34}$)\\\hline
\end{tabular}
\end{center}

Computing $T_n$ with~$n$ composite is very time consuming when~$p$ is
large, so it is important to choose the right $T_n$ quickly.  For
$p=28279$, here is a trick we used to quickly find an~$n$ such that
$\disc(T_n)$ is not divisible by~$p$.  This trick might be used to
speed up the computation for some other levels.  The key idea is to
efficiently discover which $T_n$ to compute.  Computing $T_n$ on the
full space of modular symbols is difficult, but using projections we
can compute $T_n$ on subspaces of modular symbols with small dimension
more quickly (see, e.g., \cite[\S3.5.2]{stein:phd}).  Let~$M$ be the
space of mod~$p$ modular symbols of level $p=28279$, and let
$f=\gcd(\charpoly(T_2),\deriv(\charpoly(T_2)))$.  Let~$V$ be the
kernel of $f(T_2)$ (this takes 7 minutes to compute).  If $V=0$, we
would be done, since then $\disc(T_2)\neq 0\in\F_p$.  In fact,~$V$ has
dimension~$7$.  We find the first few integers~$n$ so that the
charpoly of $T_n$ on $V$ has distinct roots, and they are $n=34$,
$47$, $53$, and $89$.  We then computed $\charpoly(T_{34})$ directly
on the whole space and found that it has distinct roots modulo~$p$.
\end{proof}

\subsection{Some Data About Weight $4$}\label{sec:k4}
The following are the valuations $d=d_4(\Gamma_0(p))$ at~$p$ of
  the discriminant of the Hecke algebras associated to
  $S_4(\Gamma_0(p))$ for $p<500$.
This data suggests a pattern, which motivates
Conjecture~\ref{conj:big} below.\vspace{1ex}
  
\hspace{-1em}\begin{minipage}[b]{0.98\textwidth}
\begin{tabular}{|l|ccccccccccccccccc|}\hline
$p$\,\, &2& 3& 5& 7& 11& 13& 17& 19& 23& 29& 31& 37& 41& 43& 47& 53& 59\\ 
$d$ &0& 0& 0& 0& 0& 2& 2& 2& 2& 4& 4& 6& 6& 6& 6& 8& 8\\\hline
$p$&61& 67& 71& 73& 79& 83& 89& 97& 101& 103& 107& 109& 113& 127& 131&  137& 139
\\
$d$ & 10& 10& 10& 12& 12& 12& 14& 16& 16& 16& 16& 18& 18& 20& 20& 22&24\\\hline
$p$ & 149& 151& 157& 163& 167& 173& 179& 181& 191& 193& 197& 199&
  211& 223& 227& 229& 233\\
$d$ & 24& 24& 26& 26& 26&28& 28& 30& 30& 32& 32& 32& 34& 36& 36& 38& 38\\ \hline
$p$ & 239& 241& 251& 257& 263& 269& 271& 277&
  281& 283& 293& 307& 311& 313& 317& 331& 337\\
$d$ & 38& 40& 40& 42& 42&44& 44& 46& 46& 46& 48& 50& 50& 52& 52& 54& 56\\\hline
$p$ & 347& 349& 353& 359& 367& 373& 379& 383& 389&397& 401& 409& 419& 421& 431& 
433& 439 \\
$d$ & 56& 58& 58& 58& 60&62& 62& 62& 65  &66& 66& 68& 68& 70& 70& 72& 72\\\hline
$p$ &  443& 449& 457& 461& 463& 467& 479& 487& 491& 499 &&&&&&&\\
$d$ &  72& 74& 76& 76& 76& 76& 78& 80& 80& 82 &&&&&&&\\\hline
\end{tabular}
\end{minipage}
\vspace{1em}

\section{Speculations}
Motivated by the promise of a pattern suggested by the table in
Section~\ref{sec:k4}, we computed $d_k(\Gamma_0(p))$ for many values
of~$k$ and~$p$. Our observations led us to the following results and
conjectures.
\begin{theorem}
  Suppose~$p$ is a prime and $k\geq 4$ is an even integer. Then
  $d_k(\Gamma_0(p))>0$ unless
\begin{align*}
  (p,k) \in \{&(2,4),(2,6),(2,8),(2,10),\\
  &(3,4),(3,6), (3,8),\\
  &(5,4), (5,6), (7,4), (11,4)\},
\end{align*}
in which case $d_k(\Gamma_0(p))=0$.
\end{theorem}

\begin{proof}
From~\cite{ribet:gerd}, mod $p$ eigenforms on $\Gamma_0(p)$ of
weight $k$ arise exactly from  mod $p$ eigenforms on $\Gamma_0(1)$ of 
weight $(k/2)(p+1)$. Moreover, there is an equality of dimensions of
vector spaces:
$$\dim S_{(k/2)(p+1)}(\Gamma_0(1)) + \dim S_{(k/2)(p+1) - (p-1)}(\Gamma_0(1))
= \dim S_{k}(\Gamma_0(p)).$$
Thus the dimension of $S_{k}(\Gamma_0(p))$ is bigger than the number
of mod $p$ eigenforms whenever 
$\dim S_{(k/2)(p+1) - (p-1)}(\Gamma_0(1))$ is non-zero. The cases of dimension
zero correspond exactly to the finite list of exceptions above,
for which one can explicitly calculate that $d_k(\Gamma_0(p)) = 0$.
\end{proof}

\medskip
Note that
for $k = 2$, however, there is a
canonical identification of spaces 
$$S_{(p+1)}(\Gamma_0(1),\Fbar_p) \simeq S_{2}(\Gamma_0(p),\Fbar_p),$$
described geometrically in~\cite{gross:duke}. For $k = 4$, the data
suggests that the discriminants $d_4(\Gamma_0(p))$ are significantly
larger than zero for large $p$, and the table above suggests a formula
of the form $2 \cdot \lfloor p/12 \rfloor$ (Not entirely
co-incidentally, this is the difference in dimension of the spaces
$S_4(\Gamma_0(p))$ and $S_{2(p+1)}(\Gamma_0(1))$).  This exact formula
is not correct, however, as evidenced by the case when $p = 139$. If
we consider the Hecke algebra $\T_4$ for $p=139$ in more detail,
however, we observe that $\T_4 \otimes \Q_{139}$ is \emph{ramified} at
$139$, and in particular contains two copies of the 
field
$\Q_{139}(\sqrt{139})$.  Just as in
the case when $k=2$ and $p=389$, there is a ``self congruence''
between the associated ramified eigenforms and their Galois
conjugates. For all other~$p$ in the range of the table, there is no
ramification, and all congruences take place between distinct
eigenforms.  Such congruences are measured by the \emph{index} of the
Hecke algebra, which is defined to be the index of $\T$ in its
normalisation $\Tw$. If we are only interested in mod $p$ congruences
(rather than mod $\ell$ congruences for $\ell \ne p$), one can
restrict to the index of $\T \otimes \Z_p$ inside its normalisation.
There is a direct relation between the discriminant and the index.
Suppose that $\T \otimes \Q_p = \prod K_i$ for certain fields
$K_i/\Q_p$ (We may assume here that $\T$ is not nilpotent, for
otherwise both the discriminant and index are infinite).  Then if
$i_p(\Gamma) = \mathrm{ord}_p([\T,\Tw])$, then
$$d_p(\Gamma) = 2 i_p(\Gamma) + \sum \mathrm{ord}_p(\Delta(K_i/\Q_p)).$$
If we now return to the example $k=4$ and $p=139$, we see
that the discrepancy from the discriminant $d_p(\Gamma_0(139)) = 24$
to the estimate
$2 \lfloor 139/12 \rfloor = 22$ is exactly accounted for by
the  two eigenforms with coefficients in $\Q_{139}(\sqrt{139})$,  which
contribute $2$ to the above formula.
 This leads us to predict that the index is exactly
given by the formula $\lfloor p/12 \rfloor$. Note that
for primes $p$ this is exactly the dimension
of $S_{p+3}(\Gamma_0(1))$.  Similar computations
lead to the following more general conjecture.


\medskip

Let~$k=2m$ be an even integer and~$p$ a prime.  Let $\T$ be the Hecke
algebra associated to $S_k(\Gamma_0(p))$ and let $\tT$ be the integral
closure of $\T$ in $\T\tensor\Q$ (which is a product of number
fields). 

\begin{conjecture}\label{conj:big}
Suppose $p\geq k-1$.  Then
  $$
  \ord_p([\tT : \T]) = \left\lfloor\frac{p}{12}\right\rfloor\cdot
  \binom{m}{2} + a(p,m),
  $$
  where
  $$
  a(p,m) =
\begin{cases}
  0 & \text{if $p\con 1\pmod{12}$,}\\
  3\cdot\ds\binom{\lceil \frac{m}{3}\rceil}{2} & \text{if $p\con 5\pmod{12}$,}\\
  2\cdot\ds\binom{\lceil \frac{m}{2}\rceil}{2} & \text{if $p\con 7\pmod{12}$,}\\
  a(5,m)+a(7,m) & \text{if $p\con 11\pmod{12}$.}
\end{cases}
$$
\end{conjecture}
Here $\binom{x}{y}$ is the binomial coefficient ``$x$ choose $y$'',
and floor and ceiling are as usual.  The conjecture is very false if
$k\gg{}p$.

When $k=2$, the conjecture specializes to the assertion that
$[\tT:\T]$ is not divisible by~$p$.  A possibly more familiar concrete
consequence of the conjecture is the following conjecture about
elliptic curves.  The modular degree of an elliptic curve~$E$ is the
smallest degree of a surjective morphism $X_0(N)\to E$, where $N$ is
the conductor of~$E$.
\begin{conjecture}
  Suppose~$E$ is an elliptic curve of prime conductor~$p$.  Then~$p$
  does not divide the modular degree $m_E$  of~$E$.
\end{conjecture}
Using the algorithm in \cite{watkins:moddeg}, M.~Watkins has computed
modular degrees of a huge number of elliptic curves of prime conductor
$p<10^7$, and not found a counterexample.  Looking at smaller data,
there is only one elliptic curve~$E$ of prime conductor $p<20000$ such
that the modular degree of~$E$ is even as big as the conductor of~$E$,
and that is a curve of conductor $13723$.  This curve has equation
$[1,1,1,-10481,408636]$, modular degree $m_E = 16176=2^4\cdot 3\cdot 337$.
The modular degree can be divisible by large primes.  For example,
there is a Neumann-Setzer elliptic curve of prime conductor $90687593$
whose modular degree is $1280092043$, which is over $14$ times as big
as $90687593$. In general, for an elliptic curve of conductor $N$,
one has the estimate $m_E \gg N^{7/6 - \epsilon}$ (see
\cite{watkins:boundmoddeg}).


\section{Conjectures Inspired by Conjecture~\ref{conj:big}}

First, some notation. Let $p$ be an odd prime.  Let $\Gamma =
\Gamma_0(p)$, and let $$S_k(R):=S_k(\Gamma)^{\new} \otimes R.$$ The
spaces $S_k$ carry an action of the Hecke algebra $\Tn_k$, and a
Fricke involution $w_p$. If $\frac{1}{2} \in R$, the space $S_k$ can
be decomposed into $+$ and $-$ eigenspaces for $w_p$.  We call the
resulting spaces $S^{+}_k$ and $S^{-}_k$ respectively.  Similarly, let
$M^+_k$ and $M^-_k$ be the $+1$ and $-1$ eigenspaces for $w_p$ on the
full spaces of new modular forms of weight~$k$ for $\Gamma_0(p)$.

It follows from \cite[Lem.~7]{atkin-lehner} (which is an explicit
formula for the trace to lower level) and the fact that $U_p$ and
$w_p$ both preserve the new subspace, that the action of the Hecke
operator $U_p$ on $S_k$ is given by the formula
$$
  U_p = -p^{(k-2)/2} w_p.
$$
This gives rise to two quotients of the Hecke algebra:
$$\T^{+} = \Tn/(U_p + p^{(k-2)/2}) \quad\text{ and }\quad \T^{-} =
\Tn/(U_p - p^{(k-2)/2}).$$ where $\T^+$ and $\T^{-}$ act on $S^{+}$
and $S^{-}$, respectively.
Recall that $\tT$ is the normalization
(integral closure) of $\T$ in $\T\tensor\Q$.  Let $\tT^{\new}$ denote
the integral closure of $\T^{\new}$ in $\T^{\new}\tensor\Q$.  
\begin{lemma} There are injections
$$\Tn \hra \T^{+}\oplus \T^{-} \hra \tTn.$$
\end{lemma}

We now begin stating some conjectures regarding the
rings $\T^{\pm}$.

\begin{conjecture} \label{conjecture:int}  Let $k < p-1$. Then
  $\T^{+}$ and $\T^{-}$ are integrally closed.  Equivalently, all
  congruences between distinct eigenforms in $S_k(\Zbar_p)$ take place
  between $+$ and $-$ eigenforms.
\end{conjecture}

Note that for $k=2$, there cannot be any congruences between $+$ and
$-$ forms because this would force $1 \equiv -1 \mod p$, which is
false, because~$p$ is odd. Thus we recover the conjecture that 
$p\nmid [\tT:\T]$ when $k=2$.  Our further conjectures go on to describe
explicitly the congruences between forms in $S_k^{+}$ and $S_k^{-}$.

\medskip

Let $E_2$ be the non-holomorphic Eisenstein series
of weight $2$. The $q$-expansion of $E_2$ is given explicitly
by
$$E_2 = 1 - 24 \sum_{n=1}^{\infty} q^n \left( \sum_{d|n} d
\right).$$
Moreover, the function $E^*_2 = E_2(\tau) - p E_2(p \tau)$ is
holomorphic of weight $2$ and level $\Gamma_0(p)$,
and moreover on $q$-expansions, $E^*_2 \equiv E_2 \mod p$.

\begin{lemma} Let $p > 3$.
Let $f \in M_k(\Gamma_0(p),\Fbar_p)$ be
a Hecke eigenform. Then $\theta f$ is an eigenform
inside $S_{k+2}(\Gamma_0(p),\Fbar_p)$.
\end{lemma}

\begin{proof} One knows that $\partial f = \theta f - k E_2 f/12$
is of weight $k+2$.  On $q$-expansions, $E_2 \equiv E^*_2 \mod p$, and
thus for $p > 3$,
$$\theta f \equiv \partial f + k E^*_2 f/12 \pmod{p}$$
is the reduction of a
weight $k+2$ form of level $p$. It is easy to see that
$\theta f$ is a cuspidal Hecke eigenform.
\end{proof}

Let us now assume Conjecture~\ref{conjecture:int} and consider the
implications for $k=4$ in more detail.  The space of modular forms
$M_2(\Gamma_0(p),\Fbar_p)$ consists precisely of $S_2$ and the
Eisenstein series $E^*_2$.  The map $\theta$ defined above induces
maps:
$$\theta: S^+_2(\Fbar_p) \rightarrow S_4(\Fbar_p),
\qquad \theta: M^{-}_2(\Fbar_p) \rightarrow S_4(\Fbar_p).$$
The images are distinct, since $\theta f = \theta g$
implies (with some care about $a_p$) that $f = g$.

\begin{conjecture} Let $f \in S_2(\Zbar_p)$ and $g \in S_4(\Zbar_p)$
be two eigenforms such that $\theta f \equiv g \mod p$. Then the
eigenvalue of $w_p$ on $f$ and $g$ have opposite signs.
\end{conjecture}

Assuming this, we get inclusions:
$$\theta S^{+}_2(\Fbar_p) \hra S^-_4(\Fbar_p),
\qquad \theta M^{-}_2(\Fbar_p) \hra S^{+}_4(\Fbar_p).$$

Now we are ready to state our main conjecture:

\begin{conjecture} There is an Hecke equivariant exact sequence:
  $$\begin{diagram} 0 &\rTo & \theta S^{+}_2(\Fbar_p) & \rTo &
    S^-_4(\Fbar_p) & \rTo & S^{+}_4(\Fbar_p) & \rTo & \theta
    M^{-}_2(\Fbar_p) & \rTo & 0. \end{diagram}$$
  Moreover, the map
  $S^{-}_4(\Fpbar) \rightarrow S^+_4(\Fpbar)$ here is the largest such
  equivariant map between these spaces.  Equivalently, a residual
  eigenform of weight $4$ and level $p$ occurs in both the $+$ and $-$
  spaces if and only if it is not in the image of $\theta$.
\end{conjecture}

Let us give some consequences of our conjectures for the index of
$\T^{\new}$ inside its normalisation.  Fix a residual representation
$\rhobar : \Gal(\Qbar/\Q) \rightarrow \GL_2(\F_q)$
 and consider the associated maximal ideal $\mathfrak{m}$
inside $\T_4$. If $\rhobar$ lies in the image of $\theta$ then our
conjecture implies that it is not congruent to any other eigenform.
If~$\rhobar$ is not in the image of $\theta$, then it should arise
exactly from a pair of eigenforms, one
inside $S^{+}_4(\Qbar_p)$ and one
inside $S^{-}_4(\Qbar_p)$. Suppose that $q = p^r$.
If there is no ramification in
$\T \otimes \Q$ over $p$ (this is often true),
 then the $+$ and $-$ eigenforms will both be
defined over the ring $W(\F_q)$ of Witt vectors of $\F_q$. Since $U_p
= p$ on $S^{-}_4$ and $-p$ on $S^{+}_4$, these forms can be at most
congruent modulo $p$. Thus the completed Hecke algebra
$(\T_4)_{\mathfrak{m}}$ is exactly
$$ \{(a,b) \in W(\F_q) \oplus W(\F_q), | a \equiv b \mod p\}.$$
One sees that this has index $q=p^r$ inside its normalisation.
Thus the (log of the)
total index is equal to $\sum r_i$ over all eigenforms
that occur inside $S^{+}_4$ and $S^{-}_4$, which from our exact
sequence we see is equal to
$$\dim S^-_4 - \dim S^+_2.$$
Conjecture~\ref{conj:big} when $k=4$,
would then follow from the equality
of dimensions:
  $$\dim S^-_4(\Fbar_p) - \dim S^+_2(\Fbar_p) = \left\lfloor
  \frac{p}{12} \right\rfloor.$$

We expect that something similar, but a little more complicated,
should happen in general.  In weight $2k$, there are mod $p^{k-r}$
congruences exactly between forms in the image of $\theta^{r-1}$ but
not of $\theta^{r}$.

\subsection{Examples}
We write small $s$'s and $m$'s for dimensions below.

Let $p=101$. Then $s^+_2 = 1$, $m^{-}_2 = 7+1 = 8$, $s^-_4 = 9$,
$s^{+}_4 = 16$. We predict the index should be $9-1 = 8 = \lfloor
101/12 \rfloor$. In the table below, we show the characteristic
polynomials of $T_2$ on $S^{-}_4$ and $S^{+}_4$, and for weight $2$,
we take the characteristic polynomial of $\theta T_2$ (or the same,
taking $F(x/2)$ where $F(x)$ is the characteristic polynomial of
$T_2$). Note that we have to add the Eisenstein series, which has
characteristic polynomial $x-1-2$, which becomes $x-6 \equiv x + 95
\mod 101$ under $\theta$.

\begin{center}
{\bf Factors of the Characteristic Polynomial of $T_2$ for $p=101$.}\vspace{1ex}
\\
\begin{tabular}{|l|l|l|l|}
\hline
$\theta S^+_2(\Fbar_{101})$ & $S^-_4(\Fbar_{101})$ & $S^+_4(\Fbar_{101})$
 & $\theta M^-_2(\Fbar_{101})$ \\
\hline
$\quad{}(x)$ &   $(x)$ & $(x+ 46)$ & $(x + 95)$ \\
& $(x + 46)$ & $(x + 95)$ & $(x^2 + 90x + 78)$ \\
& $(x^2 + 58x + 100)$ & $(x^2 + 58x + 100)$ & $(x^2 + 96x + 36)$ \\
& $(x^5 + 2x^4 + 27x^3$ & $(x^2 + 90x + 78)$ & $(x^3 + 16x^2 $ \\
& $ \ + 49x^2 + 7x + 65)$ & $(x^2 + 96x + 36)$ & $ \ + 35x + 72)$\\
 & & $(x^3 + 16x^2 + 35x + 72)$ & \\
& &  $(x^5 + 2x^4 + 27x^3  $ & \\
 & & $\ + 49x^2+ 7x + 65)$  & \\
\hline
\end{tabular}
\end{center}

Here are some further conjectures when $k > 4$.
\begin{conjecture} Let $p$ and $k$ be such that $4 < k < p-1$.
  There is an Hecke equivariant exact sequence:
$$\begin{diagram}
0 &\rTo & \theta S^{+}_{k-2}(\Fbar_p) & \rTo &
 S^-_k(\Fbar_p) & \rTo  & S^{+}_k(\Fbar_p)  & \rTo &
\theta S^{-}_{k-2}(\Fbar_p) & \rTo & 0. \end{diagram}$$
Moreover, all forms not in the image of $\theta$ contribute
maximally to the index (a factor of $p^{(k-2)/2}$). Thus
the total index should be equal to 
$$\frac{(k-2)}{2} (\dim S^+_k - \dim S^-_{k-2}) \quad + \quad\text{the index
  at level $p$ and weight $k-2$.}$$
This is the sum
$$\sum_{n=2}^{k} \frac{(2n-2)}{2} (s^+_{2n} -  s^{-}_{2n-2}).$$
\end{conjecture}

When $k = 4$, we need to add the Eisenstein series
to $S^{-}_2$ in our previous conjecture.
Note that $s^+_{k} - s^{-}_{k-2} = s^{-}_k - s^{+}_{k-2}$ for
$k > 4$ (and with $s^{-}_2$ replaced by $m^{-}_2$ when $k = 2$).
This follows from our conjectures, but can easily be proved
directly.
As an example, when $p=101$, we have $s^+_2 = 1$,
$s^-_4 = 9$, $s^+_6 = 17$, $s^-_8 = 26$, $s^+_{10} = 34$,
$s^-_{12} = 42$, $s^+_{14} = 51$, and so we would predict
the indexes $I_k$ to be as given in the following table:
\begin{center}
\begin{tabular}{|l|l|}
\hline
$k$ & $I_k?$ \\
\hline
$2$ & $0$ \\
$4$ & $8 = 8 + 0$ \\
$6$ & $24 = 24 + 0$ \\
$8$ & $51 = 48 + 3$ \\
$10$ & $83 = 80 + 3$ \\
$12$ & $123 = 120 + 3$ \\
$14$ & $177 = 168 + 9$ \\
\hline
\end{tabular}
\end{center}

This agrees with our conjectural formula, which says that
the index should be equal in this case to
$$8 \binom{k/2}{2} + 3 \binom{\lceil k/6 \rceil}{2}.$$
it also agrees with computation.


\providecommand{\bysame}{\leavevmode\hbox to3em{\hrulefill}\thinspace}
\providecommand{\MR}{\relax\ifhmode\unskip\space\fi MR }
\providecommand{\MRhref}[2]{%
  \href{http://www.ams.org/mathscinet-getitem?mr=#1}{#2}
}
\providecommand{\href}[2]{#2}

\end{document}